\documentclass[a4paper]{article}

\usepackage{amsmath,amssymb}

\newcommand{\kmin}{\underline{k}}
\newcommand{\kmax}{\overline{k}}
\newcommand{\af}{\alpha_f}
\newcommand{\ap}{\alpha_p}

\newtheorem{proposition1}{Proposition}
\newtheorem{proposition2}{Proposition}[proposition1]

\usepackage{a4wide,graphicx}

\begin{document}

\begin{center}
{\Large{\textbf{Optimized Neumann-Neumann method for the Stokes-Darcy problem}}}

\smallskip

{Marco Discacciati and Jake Robinson}

\smallskip

Department of Mathematical Sciences, Loughborough University\\ Epinal Way, Loughborough, LE11 3TU, United Kingdom\\ {m.discacciati@lboro.ac.uk}
%\and Name of Second Author \at Name, Address of Institute \email{name@email.address}
\end{center}

\medskip

\textbf{Abstract}
{A novel preconditioner of Neumann-Neumann type for the Stokes-Darcy problem is studied, where optimal weights of the local subproblems that define the preconditioner are obtained by minimizing the convergence rate of the method in the frequency space. Numerical tests show that the preconditioner is robust with respect to both the mesh size and the values of the physical parameters of the problem.}

\medskip

%%%%%%%%%
%%%%%%%%%

\section{Introduction and problem setting}\label{sect:setting}

%\vspace*{-3mm}

The Stokes-Darcy problem \cite{Discacciati:2002:MNM,Layton:2003:CFF} is a good example of multi-physics problem where splitting methods typical of domain decomposition naturally apply.
The problem is defined in a computational domain formed by a fluid region $\Omega_f$ and a porous-medium region $\Omega_p$ that are non-overlapping and separated by an interface $\Gamma$.
In $\Omega_f$, an incompressible fluid with constant viscosity and density is modelled by the dimensionless Stokes equations:
\begin{equation}\label{eq:stokes}
- \boldsymbol{\nabla} \cdot (2\mu_f \nabla^s \mathbf{u}_f - p_f \boldsymbol{I}) = \mathbf{f}_f\,, \qquad 
\nabla \cdot \mathbf{u}_f = 0\, \quad \mbox{in } \Omega_f\,,
\end{equation}
where $\mu_f = Re^{-1}$, $Re$ being the Reynolds number, $\mathbf{u}_f$ and $p_f$ are the fluid velocity and pressure, $\boldsymbol{I}$ and $\nabla^s \mathbf{u}_f = \frac{1}{2} (\nabla \mathbf{u}_f + (\nabla \mathbf{u}_f)^T)$ are the identity and the strain rate tensor, and $\mathbf{f}_f$ is an external force.
In the porous medium domain $\Omega_p$, we consider the dimensionless Darcy's model:
\begin{equation}
- \nabla \cdot (\boldsymbol{\eta}_p \nabla p_p) = f_p \qquad \mbox{in } \Omega_p\,,
\end{equation}
where $p_p$ is the fluid pressure in the porous medium, $\boldsymbol{\eta}_p$ is the permeability tensor, and $f_p$ is an external force.
The two local problems are coupled through the classical Beaver-Joseph-Saffman conditions at the interface \cite{Beavers:1967:BCN,Jager:1996:OBC,Saffman:1971:OBC}:
\begin{eqnarray}
\mathbf{u}_f \cdot \mathbf{n} = - ( \boldsymbol{\eta}_p \nabla p_p ) \cdot \mathbf{n} && \mbox{on } \Gamma, \label{eq:interf1} \\
- \mathbf{n} \cdot ( 2\mu_f \nabla^s \mathbf{u}_f - p_f \boldsymbol{I} ) \cdot \mathbf{n} = p_p && \mbox{on } \Gamma, \label{eq:interf2} \\
- ( ( 2\mu_f \nabla^s \mathbf{u}_f - p_f \boldsymbol{I} ) \cdot \mathbf{n})_\tau = \xi_f (\mathbf{u}_f)_\tau && \mbox{on } \Gamma, \label{eq:interf3}
\end{eqnarray}
where $\xi_f = \alpha_{BJ} \, (\mu_f/(\boldsymbol{\tau}\cdot \boldsymbol{\eta}_p \cdot \boldsymbol{\tau})\,)^{1/2}$, $\alpha_{BJ}$ is the Beavers-Joseph constant, $\mathbf{n}$ denotes the unit normal vector pointing outward of $\Omega_f$, while $(\mathbf{v})_{\tau}$ indicates the tangential component of any vector $\mathbf{v}$ at $\Gamma$.
Finally, we impose $\mathbf{u}_f = \mathbf{0}$ on $\Gamma_f^D$, $(2\mu_f \nabla^s\mathbf{u}_f - p_f \boldsymbol{I})\cdot\mathbf{n} = \mathbf{0}$ on $\Gamma_f^N$, $p_p = 0$ on $\Gamma_p^D$, $\mathbf{u}_p \cdot \mathbf{n}_p = 0$ on $\Gamma_p^N$, where $\Gamma_f^D \cup \Gamma_f^N = \partial\Omega_f \setminus \Gamma$ and $\Gamma_p^D \cup \Gamma_p^N = \partial\Omega_p \setminus \Gamma$.

Classical Dirichlet-Neumann type methods \cite{Quarteroni:1999:DDM} for the Stokes-Darcy problem were studied in \cite{Discacciati:2002:MNM,Discacciati:2004:CAS,Discacciati:2004:PHD} where it was pointed out that their convergence can be slow for small values of the fluid viscosity and of the porous medium permeability. Robin-Robin methods were then proposed as an alternative \cite{Discacciati:2004:PHD,Discacciati:2007:RRD,Cao:2010:FEA,Cao:2008:CSD,Chen:2011:PRR,Caiazzo:2014:CIS}, and they were analysed in the framework of optimized Schwarz methods in \cite{DGG2017,Gander:2020,DV2022}.

In this work, we focus on a Neumann-Neumann approach that allows to solve a scalar interface problem like in the case of Dirichlet-Neumann methods. This reduces the number of interface unknowns compared to the system associated with Robin-Robin iterations, and it allows to use preconditioned conjugate gradient (PCG) iterations instead of the more expensive GMRES iterations used in the Robin-Robin context (see, e.g., \cite{DGG2017}). However, to define effective Neumann-Neumann methods, the contribution of each subproblem must be suitably weighted. For single-physics problems, this is typically done using algebraic strategies that can take into account coefficient jumps across interfaces (see, e.g., \cite{Toselli:2005:DDM}). However, no clear strategies are available for multi-physics problems. In this work, we extend techniques for the analysis of optimized Schwarz methods with the aim of characterizing optimal weighting parameters to define a robust Neumann-Neumann preconditioner.

\section{Optimized Neumann-Neumann method}

Let $\alpha_f$ and $\alpha_p$ be two positive parameters: $\alpha_f, \alpha_p \in \mathbb{R}$, $\alpha_f, \alpha_p > 0$. The Neumann-Neumann method for the Stokes-Darcy problem considering the normal velocity on $\Gamma$ as interface variable reads as follows. Given $\lambda^0$ on $\Gamma$, for $m \geq 1$ until convergence,
\begin{enumerate}
\item
Find $\mathbf{u}_f^{(m)}$ and $p_f^{(m)}$ such that
\begin{equation}\label{eq:nn1}
\begin{array}{rcll}
-\boldsymbol{\nabla}\cdot(2\mu_f \nabla^s\mathbf{u}_f^{(m)} - p_f^{(m)} \boldsymbol{I}) 
= \mathbf{f}_f \, , \quad 
\nabla \cdot \mathbf{u}_f^{(m)} &=& 0 & \mbox{in } \Omega_f\,, \\
-(\mathbf{n} \cdot (2\mu_f \nabla^s \mathbf{u}_f^{(m)} - p_f^{(m)} \boldsymbol{I}))_\tau &=& \xi_f (\mathbf{u}_f^{(m)})_\tau & \mbox{on } \Gamma\,, \\
\mathbf{u}_f^{(m)}\cdot\mathbf{n} &=& \lambda^{(m-1)} & \mbox{on } \Gamma\,.
\end{array}
\end{equation}
\item
Find $p_p^{(m)}$ such that
\begin{equation}\label{eq:nn2}
\begin{array}{rcll}
-\nabla \cdot (\boldsymbol{\eta}_p \nabla p_p^{(m)}) &=& f_p & \mbox{in } \Omega_p \,,\\
-(\boldsymbol{\eta}_p \nabla p_p^{(m)}) \cdot \mathbf{n} &=& \lambda^{(m)} & \mbox{on } \Gamma\,.
\end{array}
\end{equation}
\item
Compute
\begin{equation}\label{eq:nn3}
\sigma^{(m)} = -\mathbf{n} \cdot (2\mu_f \nabla^s \mathbf{u}_f^{(m)} - p_f^{(m)} \boldsymbol{I})\cdot\mathbf{n} - p_p^{(m)} \quad \mbox{on } \Gamma\,.
\end{equation}
\item
Find $\mathbf{v}_f^{(m)}$ and $q_f^{(m)}$ such that
\begin{equation}\label{eq:nn4}
\begin{array}{rcll}
-\boldsymbol{\nabla}\cdot(2\mu_f \nabla^s\mathbf{v}_f^{(m)} - q_f^{(m)} \boldsymbol{I}) 
= \mathbf{0}  \, , \quad 
\nabla \cdot \mathbf{v}_f^{(m)} &=& 0 & \mbox{in } \Omega_f \,,\\
-(\mathbf{n} \cdot (2\mu_f \nabla^s \mathbf{v}_f^{(m)} - q_f^{(m)} \boldsymbol{I}))_\tau &=& \xi_f (\mathbf{v}_f^{(m)})_\tau & \mbox{on } \Gamma \,,\\
-\mathbf{n} \cdot (2\mu_f \nabla^s \mathbf{v}_f^{(m)} - q_f^{(m)} \boldsymbol{I})\cdot\mathbf{n} &=& \sigma^{(m)} & \mbox{on } \Gamma\,.
\end{array}
\end{equation}
\item
Find $q_p^{(m)}$ such that
\begin{equation}\label{eq:nn5}
\begin{array}{rcll}
-\nabla \cdot (\boldsymbol{\eta}_p \nabla q_p^{(m)}) &=& 0 & \mbox{in } \Omega_p \, ,\\
q_p^{(m)} &=& \sigma^{(m)} & \mbox{on } \Gamma \, .
\end{array}
\end{equation}
\item
Set
\begin{equation}\label{eq:nn6}
\lambda^{(m+1)} = \lambda^{(m)} - (\alpha_f (\mathbf{v}_f^{(m)}\cdot\mathbf{n}) + \alpha_p (\boldsymbol{\eta}_p \nabla q_p^{(m)})\cdot\mathbf{n}) \quad \mbox{on } \Gamma \, .
\end{equation}
\end{enumerate}

\noindent Problems \eqref{eq:nn1}, \eqref{eq:nn2}, \eqref{eq:nn4} and \eqref{eq:nn5} are supplemented with homogeneous boundary conditions on $\partial\Omega_f\setminus\Gamma$ and $\partial\Omega_p\setminus\Gamma$ as indicated in Sect. \ref{sect:setting}.

\subsection{Convergence analysis and optimization of the parameters}

We analyse the Neumann-Neumann method \eqref{eq:nn1}-\eqref{eq:nn6} with the aim of characterizing optimal parameters $\alpha_f$ and $\alpha_p$. To this purpose, we extend the methodology used to study optimized Schwarz methods for the Stokes-Darcy problem in \cite{DGG2017,Gander:2020,DV2022}.
Since all the problems are linear, we can study the convergence on the error equation to the zero solution when the forcing terms are $\mathbf{f}_f=\mathbf{0}$ and $f_p=0$.

We consider the simplified setting where $\Omega_f = \{(x,y)\in\mathbb{R}^2:\, x<0\}$, $\Omega_p= \{(x,y)\in\mathbb{R}^2:\, x>0\}$, $\Gamma =  \{(x,y)\in\mathbb{R}^2:\, x=0\}$, and $\mathbf{n}=(1,0)$ and $\boldsymbol{\tau}=(0,1)$. We assume $\boldsymbol{\eta}_p = \mbox{diag}(\eta_1,\eta_2)$ with constant $\eta_1\neq\eta_2$, and let $\mathbf{u}_f(x,y)= ( u_1(x,y), u_2(x,y) )^T$, $\mathbf{v}_f(x,y)= ( v_1(x,y), v_2(x,y) )^T$. In this setting, the Neumann-Neumann algorithm \eqref{eq:nn1}-\eqref{eq:nn6} becomes: given $\lambda^0$ on $\Gamma$, for $m \geq 1$ until convergence,
\begin{enumerate}
\item
Solve the Stokes problem
\begin{equation}\label{eq:nn1_ft}
\hspace*{-0.7cm}
\begin{array}{rl}
-\mu_f
\begin{pmatrix}
(\partial_{xx} + \partial_{yy}) u_1^{(m)}\\
(\partial_{xx} + \partial_{yy}) u_2^{(m)}
\end{pmatrix}
+
\begin{pmatrix}
\partial_x p_f^{(m)} \\
\partial_y p_f^{(m)}
\end{pmatrix}
=0, \;
\partial_x u_1^{(m)} + \partial_y u_2^{(m)} = 0,
 & \;\mbox{in } (-\infty,0) \times \mathbb{R}\,, \\[10pt]
%\partial_x u_1^{(m)} + \partial_y u_2^{(m)} &=& 0 & \mbox{in } (-\infty,0) \times \mathbb{R} \,,\\
-\mu_f (\partial_x u_2^{(m)} + \partial_y u_1^{(m)}) = \xi_f \, u_2^{(m)}, \quad u_1^{(m)} = \lambda^{(m)}, & \;\mbox{on } \{ 0 \} \times \mathbb{R} \,.\\
%u_1^{(m)} &=& \lambda^{(m)} & \mbox{on } \{ 0 \} \times \mathbb{R} \,.
\end{array}
\end{equation}
\item
Solve Darcy's problem
\begin{equation}\label{eq:nn2_ft}
\begin{array}{rcll}
-( \eta_1 \, \partial_{xx} + \eta_2\, \partial_{yy} )\, p_p^{(m)} 
&=& 0 & \mbox{in } (0,+\infty) \times \mathbb{R}\,, \\
-\eta_1 \partial_x p_p^{(m)} &=& \lambda^{(m)} & \mbox{on } \{ 0 \} \times \mathbb{R} \,.
\end{array}
\end{equation}

\item
Compute
\begin{equation}\label{eq:nn3_ft}
\sigma^{(m)} = -2\mu_f \partial_x u_1^{(m)} + p_f^{(m)} - p_p^{(m)} \quad \mbox{on } \{ 0 \} \times \mathbb{R}\,.
\end{equation}
\item
Solve the Stokes problem
\begin{equation}\label{eq:nn4_ft}
\hspace*{-9mm}
\begin{array}{rl}
-\mu_f
\begin{pmatrix}
(\partial_{xx} + \partial_{yy}) v_1^{(m)}\\
(\partial_{xx} + \partial_{yy}) v_2^{(m)}
\end{pmatrix}
+
\begin{pmatrix}
\partial_x q_f^{(m)} \\
\partial_y q_f^{(m)}
\end{pmatrix}
=0 ,\; %& \mbox{in } (-\infty,0) \times \mathbb{R} \,,\\
\partial_x v_1^{(m)} + \partial_y v_2^{(m)} = 0\,, & \; \mbox{in } (-\infty,0) \times \mathbb{R} \,,\\
-\mu_f (\partial_x v_2^{(m)} + \partial_y v_1^{(m)}) = \xi_f \, v_2^{(m)}\,, & \; \mbox{on } \{ 0 \} \times \mathbb{R} \,,\\
-2\mu_f \partial_x v_1^{(m)} + q_f^{(m)} = \sigma^{(m)}\,, & \; \mbox{on } \{ 0 \} \times \mathbb{R} \,.
\end{array}
\end{equation}
\item
Solve Darcy's problem
\begin{equation}\label{eq:nn5_ft}
\begin{array}{rcll}
-( \eta_1 \, \partial_{xx} + \eta_2\, \partial_{yy} )\, q_p^{(m)} 
&=& 0 & \mbox{in } (0,+\infty) \times \mathbb{R}\,, \\
q_p^{(m)} &=& \sigma^{(m)} & \mbox{on } \{ 0 \} \times \mathbb{R}\,.
\end{array}
\end{equation}
\item
Set
\begin{equation}\label{eq:nn6_ft}
\lambda^{(m+1)} = \lambda^{(m)} - (\af v_1^{(m)} + \ap \, \eta_1 \partial_x q_p^{(m)}) \quad \mbox{on } \{ 0 \} \times \mathbb{R} \, .
\end{equation}
\end{enumerate}

For the convergence analysis, we consider the Fourier transform in the direction tangential to the interface (corresponding to the $y$ variable): 
\[
\mathcal{F} : w(x,y) \mapsto \widehat{w}(x,k) = \int_{\mathbb R}
e^{-iky}w(x,y) \,dy\,, \qquad \forall w(x,y)\in  L^2(\mathbb{R}^2)\,,
\]
where $k$ is the frequency variable.  We quantify the error in the frequency space between two successive approximations $\widehat{\lambda}^{m+1}$ and $\widehat{\lambda}^{m}$ at $\Gamma$ and characterize the reduction factor at iteration $m$ for each frequency $k$. Finally, we identify optimal values of $\af$ and $\ap$ by minimizing the reduction factor at each iteration over all the relevant Fourier modes.

\begin{proposition1}
Let $\eta_p = \sqrt{\eta_1\eta_2}$. The reduction factor of algorithm \eqref{eq:nn1_ft}-\eqref{eq:nn5_ft} does not depend on the iteration $m$, and it is given by $|\rho(\alpha_f,\alpha_p,k)|$ with
\begin{equation}\label{eq:rhog}
\rho(\alpha_f,\alpha_p,k) = 
1-\alpha_p (1 + 2\,\mu_f \,\eta_p\, k^2) - \alpha_f ( 1 + ( 2\,\mu_f\,\eta_p\,k^2 \,)^{-1} ) \, .
\end{equation} 
\end{proposition1}

%\begin{proof}
\emph{Proof.} Following the same steps of the proof of Proposition 3.1 of \cite{DGG2017}, we find
\begin{equation*}
%\begin{array}{c}
%\displaystyle
 \widehat{u}_1^{(m)} (x,k) = \left( U_1^{(m)}(k)  + \frac{P^{(m)}(k)}{2\mu_f}\,x \right) e^{|k|x}\, , %\\
%\displaystyle
%\widehat{p}_f^{(m)}(x,k) = P^{(m)}(k) \,e^{|k|x} \,,
\qquad
\widehat{p}_p^{(m)}(x,k) = \Phi^{(m)}(k) \,e^{- \sqrt{\frac{\eta_2}{\eta_1}}|k|x}\,,
%\end{array}
\end{equation*}
and $\widehat{p}_f^{(m)}(x,k) = P^{(m)}(k) \,e^{|k|x}$. The interface conditions \eqref{eq:nn1_ft}$_4$ and \eqref{eq:nn2_ft}$_2$ give $U_1^{(m)}(k) = \widehat{\lambda}^{(m)}$ and $\Phi^{(m)}(k) = \frac{\widehat{\lambda}^{(m)}}{\eta_p|k|}$.
Then, using the Fourier transform of \eqref{eq:nn3_ft}, we can obtain
\[
  \widehat{\sigma}^{(m)} = - ( 2\mu_f |k| + (\eta_p|k|)^{-1} ) \, \widehat{\lambda}^{(m)} \, .
\]
Proceeding in analogous way, the solutions of problems \eqref{eq:nn4_ft} and \eqref{eq:nn5_ft} become
\begin{equation*}
%\begin{array}{c}
%\displaystyle
 \widehat{v}_1^{(m)} (x,k) = \left( \overline{P}^{(m)}(k)\, x -\frac{\widehat{\sigma}^{(m)}}{|k|} \right) \frac{e^{|k|x}}{2\mu_f}\,, %\\
%\displaystyle
\qquad
\widehat{q}_p^{(m)} (x,k) = \widehat{\sigma}^m e^{-\sqrt{\frac{\eta_2}{\eta_1}}|k|x} \,,
\end{equation*}
and $\widehat{q}_f^{(m)} (x,k) = \overline{P}^m(k) e^{|k|x}$.
%\end{array}
%
Substituting into the Fourier transform of \eqref{eq:nn6_ft}, we find $\widehat{\lambda}^{(m+1)} = \rho(\alpha_f,\alpha_p,k) \widehat{\lambda}^{(m)}$ with $\rho(\alpha_f,\alpha_p,k)$ defined in \eqref{eq:rhog}.
$\Box$
%\end{proof}

Using a classical approach in optimized Schwarz methods, we now aim at optimizing the parameters $\alpha_f$ and $\alpha_p$ by minimizing the reduction factor for all the relevant frequencies $k$ with $0<\kmin \leq |k| \leq \kmax$, where $\kmin$ and $\kmax$ are the minimum and maximum relevant frequencies, respectively, with $\kmin = \pi/L$ ($L$ being the length of the interface) and $\kmax = \pi/h$ ($h$ being the size of the mesh). Since the function $\rho(\alpha_f,\alpha_p,k)$ is even with respect to $k$, we only consider $k>0$ without loss of generality, and we proceed to solve the min-max problem
\begin{equation}\label{eq:minmax}
  \min_{\alpha_f,\alpha_p>0} \max_{k \in [\kmin,\kmax]} \, |\rho(\alpha_f,\alpha_p,k)| \, .
\end{equation}
The following result holds.

\begin{proposition2}
The solution of the min-max problem \eqref{eq:minmax} is given by
\begin{equation}\label{eq:afap}
\begin{array}{l}
\displaystyle
\alpha_f^{NN} = (2\,\mu_f\,\eta_p\,\kmin\,\kmax)^2 \, ( \, 1 + (2\,\mu_f\,\eta_p\,\kmin\,\kmax)^2 + \mu_f \, \eta_p \, (\kmin + \kmax)^2 \, )^{-1}  \, ,\\[12pt]
\displaystyle
%\alpha_p^* = \frac{1}{1 + (2\,\mu_f\,\eta_p\,\kmin\,\kmax)^2 + \mu_f \, \eta_p \, (\kmin + \kmax)^2} \, .
\alpha_p^{NN} = (\, 1 + (2\,\mu_f\,\eta_p\,\kmin\,\kmax)^2 + \mu_f \, \eta_p \, (\kmin + \kmax)^2 \, )^{-1} \, .
\end{array}
\end{equation}
Moreover, $|\rho(\alpha_f^{NN},\alpha_p^{NN},k)|<1$ for all $k \in [\kmin,\kmax]$, and, asymptotically, when $h\to 0$,
\begin{eqnarray*}
%\alpha_f^{NN} &=& 4\pi^2\mu_f\eta_p \, C_{NN} - 8\pi^2\mu_f \eta_p \, L \, C_{NN}^{2}\, h + O(h^2) \\
\alpha_f^{NN} &=& 4\pi^2\mu_f\eta_p \, C_{NN} \, (1 - 2\,L \,C_{NN}\, h) + O(h^2) \\
\alpha_p^{NN} &=& L^2(\pi^2\mu_f\eta_p)^{-1}C_{NN}\, h^2 + O(h^3) \\
\rho(\alpha_f^{NN},\alpha_p^{NN},\kmax) &=& -L^2\, C_{NN} + (8\pi^2\mu_f\eta_pL + 4L^3)\,C_{NN}^2\, h + O(h^2)\,,
\end{eqnarray*}
with $C_{NN} = (4\pi^2\mu_f\,\eta_p+L^2)^{-1}$.
\end{proposition2}

%\begin{proof}
\emph{Proof.}
For all $\alpha_f,\alpha_p>0$, $\lim_{k\to 0} \rho(\alpha_f,\alpha_p,k) = \lim_{k\to \infty} \rho(\alpha_f,\alpha_p,k) = -\infty$, and the function $\rho(\alpha_f,\alpha_p,k)$ has a local maximum at
$k^* = ( \alpha_f / (\alpha_p \, (2\,\mu_f\,\eta_p)^2 )\,)^{1/4}$
%\[
%  k^* = \left( \frac{\alpha_f}{\alpha_p} \, \frac{1}{(2\,\mu_f\,\eta_p)^2} \right)^{1/4}
%\]
where
\begin{equation}\label{eq:ks}
  \rho(\alpha_f,\alpha_p,k^*) = 1 - (\sqrt{\af}+\sqrt{\ap})^2\,.
\end{equation}
We distinguish two cases.

\smallskip

\noindent \emph{Case 1:} $\sqrt{\af}+\sqrt{\ap} \geq 1$. In this case, $\rho(\alpha_f,\alpha_p,k) \leq 0$ for all $\kmin \leq k \leq \kmax$, and $\rho(\alpha_f,\alpha_p,k) = 0$ if $\sqrt{\af}+\sqrt{\ap} = 1$. Taking $\sqrt{\af}+\sqrt{\ap} = 1$ would result in a null convergence rate for $k=k^*$, and we could then choose $\af$ and $\ap$ by imposing $|\rho(\af,\ap,\kmin)| = |\rho(\af,\ap,\kmax)|$ (which would also ensure that $\kmin < k^* < \kmax$). This approach leads to $\ap = (1+2\,\mu_f\,\eta_p\,\kmin\,\kmax)^{-2}$ and $\af = (2\,\mu_f\,\eta_p\,\kmin\,\kmax)^{2}(1+2\,\mu_f\,\eta_p\,\kmin\,\kmax)^{-2}$, but, unfortunately, it does not guarantee that $|\rho(\af,\ap,k)|<1$ for all $k \in [\kmin,\kmax]$, which would be true when $1+2\,\mu_f\,\eta_p \, \kmin \, \kmax > \sqrt{2\, \mu_f \, \eta_p} \, (\kmax - \kmin)$.

\smallskip

\noindent \emph{Case 2:} $0 < \sqrt{\af}+\sqrt{\ap} < 1$. In this case, $\rho(\af,\ap,k^*) >0$, and the function $\rho(\af,\ap,k)$ has two positive zeros

\smallskip

\centerline{%
  $k_{1,2} = ( 1-\af-\ap \pm ( \, (1 - \af - \ap)^2 - 4\, \af\,\ap)^{1/2})^{1/2} \, / \, (4\, \mu_f \, \eta_p \, \ap)^{1/2}$ \, ,}
  
\smallskip
  
\noindent whose position depends on the values of $\af$ and $\ap$. Therefore, we proceed by equioscillation and we look for $\af$ and $\ap$ such that $-\rho(\af,\ap,\kmin) = \rho(\af,\ap,k^*)$ and $-\rho(\af,\ap,\kmax) = \rho(\af,\ap,k^*)$. This gives the values \eqref{eq:afap}. Simple algebraic manipulations permit to verify that, for such values of the parameters, $k^* = (\kmin\,\kmax)^{1/2}$ so that $\kmin < k_1 < k^* < k_2 < \kmax$. Moreover, $|\rho(\af,\ap,k)| \leq \rho(\af,\ap,k^*)$ for all $\kmin \leq k \leq \kmax$ and, owing to \eqref{eq:ks}, we can conclude that $|\rho(\af,\ap,k)| < 1$ for all frequencies of interest.
%\end{proof}
$\Box$

%\vspace*{-3mm}

\section{Numerical results}

%\vspace*{-3mm}

We consider a finite element approximation based on the inf-sup stable $\mathbb{Q}_2-\mathbb{Q}_1$ Taylor-Hood elements \cite{Boffi:2013:MFE} for Stokes, and $\mathbb{Q}_2$ elements Darcy. 
Denoting by the indices $I_f$, $I_p$ and $\Gamma$ the degrees of freedom in $\Omega_f$, $\Omega_p$ and on $\Gamma$, respectively, the algebraic form of the discrete Stokes-Darcy problem \eqref{eq:stokes}--\eqref{eq:interf3} becomes
\begin{equation}\label{eq:system}
  \begin{pmatrix}
  A^f_{I_fI_f} & A^f_{I_f\Gamma} & G^f_{I_f} & 0 & 0 \\[3pt]
  A^f_{\Gamma I_f} & A^f_{\Gamma\Gamma} & G^f_{\Gamma} & 0 & C_{fp} \\[3pt]
  (G^f_{I_f})^T & (G^f_\Gamma)^T & 0 & 0 & 0 \\
  0 & 0 & 0 & A^p_{I_pI_p} & A^p_{I_p\Gamma} \\[3pt]
  0 & -C_{fp}^T & 0 & A^p_{\Gamma I_p} & A^p_{\Gamma\Gamma}
  \end{pmatrix}
  \begin{pmatrix}
  \mathbf{u}_{f,I_f} \\
  \mathbf{u}_{f,\Gamma} \\
  \mathbf{p}_f \\
  \mathbf{p}_{p,I_p} \\
  \mathbf{p}_{p,\Gamma}  
  \end{pmatrix}
  =
  \begin{pmatrix}
  \mathbf{f}_{f,I_f} \\
  \mathbf{f}_{f,\Gamma} \\
  \mathbf{0} \\
  \mathbf{f}_{p,I_p} \\
  \mathbf{f}_{p,\Gamma}
  \end{pmatrix}\,,
\end{equation}
where $\mathbf{u}_{f,\Gamma}$ denotes the vector of degrees of freedom of the normal velocity on $\Gamma$. 
The Schur complement system with respect to $\mathbf{u}_{f,\Gamma}$ is
\begin{equation}\label{eq:schur}
  (\Sigma_f + \Sigma_p) \, \mathbf{u}_{f,\Gamma} = \mathbf{b}_{\Gamma}
\end{equation}
where $\Sigma_f$ and $\Sigma_p$ are the symmetric and positive definite matrices (see \cite{Discacciati:2004:PHD}):
\begin{equation*}
\Sigma_f = A^f_{\Gamma\Gamma} 
               - \begin{pmatrix}
                    A^f_{\Gamma I_f} &  G^f_{\Gamma} 
                 \end{pmatrix}
                 \begin{pmatrix}
                    A^f_{I_f I_f} &  G^f_{I_f} \\[3pt]
                    (G^f_{I_f})^T &  {0}
                 \end{pmatrix}^{-1}
                 \begin{pmatrix}
                    A^f_{I_f\Gamma} \\[3pt]
                    (G^f_{\Gamma})^T
                 \end{pmatrix}\, ,
\end{equation*}
\begin{equation*}
\Sigma_p = \begin{pmatrix} 
                    0 & {C}_{fp}
               \end{pmatrix} 
               \begin{pmatrix}
                   A^p_{I_p I_p}    & A^p_{I_p \Gamma} \\[3pt]
                   A^p_{\Gamma I_p} & A^p_{\Gamma \Gamma}
               \end{pmatrix}^{-1}
               \begin{pmatrix} 
                    0 \\ 
                    {C}_{fp}^T
               \end{pmatrix} \, .
\end{equation*}

Following a classical approach in domain decomposition (see, e.g., \cite{Quarteroni:1999:DDM,Discacciati:2004:PHD}), the Neumann-Neumann method \eqref{eq:nn1}--\eqref{eq:nn5} can be equivalently reformulated as a Richardson method for the Schur complement system \eqref{eq:schur} with preconditioner
\begin{equation}\label{eq:PNN}
P = \af \, \Sigma_f^{-1} + \ap\, \Sigma_p^{-1}\,.
\end{equation}
The PCG method with preconditioner $P$ can then be used to solve \eqref{eq:schur}.

\smallskip

We consider the computational domains $\Omega_f = (0,0.5) \times (1,1.5)$ and $\Omega_p = (0,0.5) \times (0.5,1)$ so that $\Gamma = (0,0.5) \times \{1\}$, and we choose the forces $\mathbf{f}_f$ and $f_p$ and the boundary conditions in such a way that the Stokes-Darcy problem has analytic solution $\mathbf{u}_f = (\sqrt{\eta_p}, \alpha_{BJ}x)^T$, $p_f = 2\,\mu_f\,(x+y-1) + (3\,\eta_p)^{-1}$, and $p_p = \eta_p^{-1} (-\alpha_{BJ} x (y-1) + y^3/3-y^2+y) + 2\,\mu_fx$. The computational meshes are structured and characterized by $h=0.1 \times 2^{1-j}$, $j=1,\ldots,4$, with 11, 21, 41, and 81 interface unknowns, respectively. We consider four configurations of physically significant dimensionless problem parameters (see also \cite{DV2022}): \emph{(a)} $\mu_f = 10$, $\eta_p = 4 \times 10^{-10}$; \emph{(b)} $\mu_f = 1$, $\eta_p = 4 \times 10^{-7}$; \emph{(c)} $\mu_f = 10$, $\eta_p = 4 \times 10^{-9}$; \emph{(d)} $\mu_f = 0.2$, $\eta_p = 2 \times 10^{-7}$.

Table \ref{table:results1} reports the computed values of the optimal parameters $\af^{NN}$ and $\ap^{NN}$ \eqref{eq:afap} and the number of CG iterations with preconditioner \eqref{eq:PNN} and without preconditioner (in brackets). For comparison, we indicate also the values of the optimal parameters $\af^{RR}$ and $\ap^{RR}$ and the number of GMRES iterations obtained with the optimized Schwarz (Robin-Robin) method studied in \cite{DGG2017}. (Notice that $\af^{NN} \sim c_f^0 + c_f^1 h$ and $\ap^{NN} \sim c_p^0 + c_p^1 h$ when $h \to 0$ for suitable constants $c_f^0$, $c_f^1$, $c_p^0$ and $c_p^1$ that depend on $\mu_f$, $\eta_p$ and $L$.)

\begin{table}[bht]
\caption{Optimal parameters $\af^{NN}$ and $\ap^{NN}$ and number of PCG iterations, and optimal parameters $\af^{RR}$ and $\ap^{RR}$ for the Robin-Robin method with corresponding GMRES iterations (\emph{tol}$\;=10^{-9}$).}\label{table:results1}
\begin{center}
\resizebox{\columnwidth}{!}{%
\begin{tabular}{ccllcccllc}
Case \;& Mesh & \multicolumn{1}{c}{$\af^{NN}$} & \multicolumn{1}{c}{$\ap^{NN}$} &\phantom{a}& \multicolumn{2}{c}{PCG iter} & \multicolumn{1}{c}{$\af^{RR}$} & \multicolumn{1}{c}{$\ap^{RR}$} & GMRES iter \\
\hline
(a) & $h_1$ & $9.97 \times 10^{-12}$\; & $1.00 \times 10^{+0}$ & & $2$ & $(12)$ & $7.23 \times 10^{+7}$\; & $6.91 \times 10^{+2}$\; & $4$ \\
    & $h_2$ &   $3.99 \times 10^{-11}$   &   $1.00 \times 10^{+0}$ & & $2$ & $(17)$ &   $3.79 \times 10^{+7}$   &   $1.32 \times 10^{+3}$   & $4$ \\
    & $h_3$ &   $1.60 \times 10^{-10}$   &   $1.00 \times 10^{+0}$ & & $3$ & $(22)$ &   $1.94 \times 10^{+7}$   &   $2.58 \times 10^{+3}$   & $4$ \\
    & $h_4$ &   $6.38 \times 10^{-10}$   &   $9.99 \times 10^{-1}$ & & $3$ & $(31)$ &   $9.83 \times 10^{+6}$   &   $5.09 \times 10^{+3}$   & $4$ \\
\hline
(b) & $h_1$ & $9.96 \times 10^{-8}$ & $9.98 \times 10^{-1}$ & & $3$ & $(12)$ & $7.24 \times 10^{+4}$ & $6.91 \times 10^{+1}$ & $6$ \\
    & $h_2$ & $3.96 \times 10^{-7}$ & $9.93 \times 10^{-1}$ & & $4$ & $(17)$ & $3.80 \times 10^{+4}$ & $1.32 \times 10^{+2}$ & $6$ \\
    & $h_3$ & $1.55 \times 10^{-6}$ & $9.74 \times 10^{-1}$ & & $4$ & $(24)$ & $1.96 \times 10^{+4}$ & $2.55 \times 10^{+2}$ & $8$ \\
    & $h_4$ & $5.78 \times 10^{-6}$ & $9.06 \times 10^{-1}$ & & $5$ & $(30)$ & $1.03 \times 10^{+4}$ & $4.86 \times 10^{+2}$ & $8$ \\
\hline
(c) & $h_1$ & $9.97 \times 10^{-10}$ & $1.00 \times 10^{+0}$ & & $3$ & $(12)$ & $7.23 \times 10^{+6}$ & $6.91 \times 10^{+2}$ & $4$ \\
    & $h_2$ & $3.99 \times 10^{-9}$  & $9.99 \times 10^{-1}$ & & $3$ & $(17)$ & $3.79 \times 10^{+6}$ & $1.32 \times 10^{+3}$ & $4$ \\
    & $h_3$ & $1.59 \times 10^{-8}$  & $9.97 \times 10^{-1}$ & & $3$ & $(24)$ & $1.94 \times 10^{+6}$ & $2.57 \times 10^{+3}$ & $6$ \\
    & $h_4$ & $6.32 \times 10^{-8}$  & $9.90 \times 10^{-1}$ & & $4$ & $(30)$ & $9.87 \times 10^{+5}$ & $5.06 \times 10^{+3}$ & $6$ \\
\hline
(d) & $h_1$ & $2.49 \times 10^{-10}$ & $1.00 \times 10^{+0}$ & & $2$ & $(12)$ & $7.23 \times 10^{+5}$ & $3.46 \times 10^{+1}$ & $4$ \\
    & $h_2$ & $9.97 \times 10^{-10}$ & $1.00 \times 10^{+0}$ & & $3$ & $(17)$ & $3.79 \times 10^{+5}$ & $6.60 \times 10^{+1}$ & $4$ \\
    & $h_3$ & $3.98 \times 10^{-9}$  & $9.99 \times 10^{-1}$ & & $3$ & $(22)$ & $1.94 \times 10^{+5}$ & $1.29 \times 10^{+2}$ & $6$ \\
    & $h_4$ & $1.59 \times 10^{-8}$  & $9.95 \times 10^{-1}$ & & $4$ & $(29)$ & $9.85 \times 10^{+4}$ & $2.54 \times 10^{+2}$ & $6$ \\
\hline
\end{tabular}
}
\end{center}
\end{table}

The number of PCG iterations using optimized parameters $\af^{NN}$, $\ap^{NN}$ is almost independent of both the mesh size and of the values of $\mu_f$ and $\eta_p$.

Moreover, the optimized Neumann-Neumann method performs better than the Robin-Robin method with lower computational cost per iteration. We also observe that, considering the Robin interface conditions (3.3)$_4$ and (3.4)$_4$ in \cite{DGG2017} and the values of $\af^{RR}$ and $\ap^{RR}$ (especially, the large values of $\af^{RR}$), the Robin-Robin method actually behaves like a Dirichlet-Robin method with interface condition on the normal velocity $\mathbf{u}_f \cdot \mathbf{n}$ for the Stokes problem. This confirms that condition \eqref{eq:nn1}$_3$ in the Neumann-Neumann algorithm is a valid choice for the Stokes problem.

Finally, the optimal values $\af^{NN}$, $\ap^{NN}$ suggest that the preconditioner \eqref{eq:PNN} behaves like $P \approx \Sigma_p^{-1}$. Thus, while $\Sigma_f^{-1}$ is an effective preconditioner for large values of $\mu_f$ and $\eta_p$ (see \cite{Discacciati:2004:PHD,Discacciati:2004:CAS}), $\Sigma_p^{-1}$ is a much better choice for small values, which is the case in most applications. This can lead to a Dirichlet-Neumann-type method different from the one in \cite{Discacciati:2004:PHD,Discacciati:2004:CAS} that will be discussed in a future work.

%\section{Conclusions}
%
%Using the approach to analyse optimized Schwarz methods, we characterized an optimal preconditioner of Neumann-Neumann type for the Schur complement system associated with the coupled Stokes-Darcy problem when the normal velocity is the interface variable. The preconditioner is robust with respect to both the mesh size and the values of the physical parameters of the problem.

%\vspace*{-2mm}

\bigskip

\emph{Acknowledgement.}
The first author acknowledges funding by the EPSRC grant EP/V027603/1.
%\end{acknowledgement}

%\vspace*{-7mm}

\bibliographystyle{plain}
\bibliography{DR}
\end{document}